\documentclass[12pt, reqno]{article}
\usepackage{amssymb}
\usepackage{amsmath}
\usepackage{amscd}
\usepackage{delarray}
\usepackage{verbatim}
\usepackage{bm}

\makeatletter

       \@addtoreset{equation}{section}
\makeatother       
\begin{document}
\bibliographystyle{plain}
\newtheorem{main}{Main theorem}
\renewcommand{\themain}{}
\newtheorem{thm}{Theorem}[section]
\newtheorem{prop}[thm]{Proposition}
\newtheorem{lemm}[thm]{Lemma}
\newtheorem{defn}[thm]{Definition}
\newtheorem{cor}[thm]{Corollary}
\newtheorem{rem}[thm]{Remark}
\newtheorem{exmp}[thm]{Example}
\newtheorem{ques}[thm]{Question}
\newtheorem{prob}[thm]{Problem}
\def\Vec#1{\mbox{\boldmath $#1$}}
\newfont{\bg}{cmr10 scaled\magstep4}
\newcommand{\bigzerol}{\smash{\hbox{\bg 0}}}
\newcommand{\bigzerou}{\smash{\lowerl.7ex\hbox{\bg 0}}} 
\renewcommand{\thefootnote}{\fnsymbol{footnote}}
\title{Linear continuous surjections of $C_{p}$-spaces over 
compacta}
\author{Kazuhiro Kawamura and Arkady Leiderman}
\date{}

\maketitle
\begin{abstract}
Let $X$ and $Y$ be compact Hausdorff spaces and suppose that 
there exists a linear continuous surjection $T:C_{p}(X) \to C_{p}(Y)$, 
where $C_{p}(X)$ denotes the space of all real-valued continuous 
functions on $X$ endowed with the pointwise convergence topology.  We 
prove that $\dim X=0$ implies $\dim Y = 0$. This generalizes 
a previous theorem \cite[Theorem 3.4]{LLP} for compact metrizable 
spaces.  Also we point out that the function space $C_{p}(P)$ over 
the pseudo-arc $P$ admits no densely defined linear continuous operator 
$C_{p}(P) \to C_{p}([0,1])$ with a dense image. 
\end{abstract}

\footnotetext{
The authors are supported by JSPS KAKENHI Grant Number 26400080.
The visit of the second author to University of Tsukuba in August, 
2015 was supported by the grant above.
}

\section{Introduction and Results}

For a Tychonoff space $X$, $C_{p}(X)$ denotes the space of all 
continuous real-valued functions on $X$ endowed with the pointwise convergence topology. The relationship between the topology of $X$ 
and linear topological properties of 
$C_{p}(X)$ is a subject of extensive research.  A theorem of Pestov \cite{Pestov} plays 
the fundamental role in this study:
if $C_{p}(X)$ and $C_{p}(Y)$ are linearly homeomorphic for Tychonoff 
spaces $X$ and $Y$, then we have 
the equality $\dim X = \dim Y$.  The theorem was first proved by 
Pavlovski{\v i} for compact metrizable spaces \cite{Pavlovskii}.
A natural 
question arises whether the inequality $\dim Y \leq \dim X$ holds for Tychonoff spaces $X$ and $Y$ whenever 
there exists a linear continuous surjection $C_{p}(X) \to C_{p}(Y)$ 
\cite[Problem 1046,1047]{ArhangProb}.
This was answered negatively in \cite{LMP},\cite{LLP} even for compact metrizable spaces and the results were recently refined in \cite{Levin}.  The exception is the zero-dimensional case 
\cite[Theorem 3.4]{LLP}: if 
there exists a linear continuous surjection $C_{p}(X) \to C_{p}(Y)$ 
for compact metrizable spaces $X$ and $Y$, then $\dim X = 0$ implies 
$\dim Y = 0$. The present paper extends the above theorem to all compact Hausdorff spaces.

\begin{thm}
Let $X$ and $Y$ be compact Hausdorff spaces and suppose that 
there exists a linear continuous surjection $T:C_{p}(X) \to C_{p}(Y)$.
If $\dim X = 0$, then we have $\dim Y = 0$.
\end{thm}

Our proof is based on the spectral theorem of Shchepin \cite{Scepin} \cite{Chigogidze} which allows us 
to reduce our consideration to that on compact metrizable spaces (Proposition 2.2). 
Then a slight modification of \cite[Theorem 3.4]{LLP} supplies the desired result (Proposition 2.1).

The proof of Proposition 2.1 is applied to obtain more information 
on the existence of linear continuous surjections.  It is known that 
for each finite- 
dimensional compact metrizable space $X$, there exists a linear continuous 
surjection $C_{p}([0,1]) \to C_{p}(X)$ \cite{LMP} and the map may 
be constructed to be open \cite{Levin}.  The assumption of finite-  
dimensionality cannot be dropped since there exists no linear continuous surjection 
$C_{p}([0,1]) \to C_{p}([0,1]^{\omega})$ \cite[Remark 4.6]{LMP}.  
We see from the next proposition that 
$[0,1]$ in these results cannot be replaced 
by the pseudo-arc $P$, the topologically unique 
hereditarily indecomposable 
continuum which is the limit of an inverse sequence of $[0,1]$ (see 
a survey article \cite{Lewis}).  
Here a continuum means a compact connected metrizable space and 
a continuum $X$ is said to be {\it hereditarily indecomposable} 
if each subcontinuum is not the union of two proper subcontinua. 
A compact metrizable 
space is called a {\it Bing compactum} if each connected component is either 
herediarily indecomposable or is a singleton.
A compact metrizable space $Y$ is said to be {\it hereditarily locally connected} 
if each subcontinuum of $Y$ is locally connected.  Examples are  
$[0,1]$ and more generally dendrites (\cite{Nadler}).

\begin{prop}
Let $X$ be a Bing compactum, $Y$ be a hereditarily locally connected 
compact metrizable space and let $T:C_{p}(X) \to C_{p}(Y)$ 
be a densely defined linear continuous operator with a dense image.
Then we have $\dim Y = 0$.  
In particular there exists no linear continuous surjection 
$C_{p}(P) \to C_{p}([0,1])$.
\end{prop}

\vspace{5mm}
For a linear operator $T:E \to F$, $D(T)$ and 
$R(T)$ denote the domain of $T$ and the image of $T$ respectively.  
An operator $T:E \to F$ is said to be {\it densely defined} if 
$D(T)$ is dense in $E$.  For a continuous map $\varphi:X \to Y$, 
the induced operator $\varphi^{\sharp}:C_{p}(Y) \to C_{p}(X)$ is 
defined by
\[
\varphi^{\sharp}(f) = f \circ \varphi,~~ f \in C_{p}(Y).
\] 

\section{Proofs}

\begin{prop}
Let $X$ and $Y$ be compact metrizable spaces and let 
$T:E \to F$ be a linear continous surjection defined on 
a dense subspace $E$ of $C_{p}(X)$ onto a dense subspace 
$F$ of $C_{p}(Y)$.  If $\dim X = 0$, then we have $\dim Y = 0$.
\end{prop} 

The compactness 
assumption cannot be dropped in the above theorem as is demonstrated 
by the following example. 

\begin{exmp}
There exists a densely defined surjective linear operator 
$T:C_{p}(C') \to C_{p}([0,1])$, where 
$C'$ is a $G_{\delta}$ subset of the Cantor set $C$ and hence 
is a Polish space.
\end{exmp}

Let $\varphi:C \to [0,1]$ be a continuous 2-to-1 map of the Cantor set 
$C$ onto $[0,1]$ such that the set 
$\{ y \in [0,1]~|~ |\varphi^{-1}(y)| = 2 \}$ is at most countable.
We then obtain a subset $C'$ of $C$, by removing a countable subset 
from $C$, such that
$\varphi '= \varphi|C' \to [0,1]$ is a continuous bijection.
The map induces a dense embedding
$(\varphi')^{\sharp}:C_{p}([0,1]) \to C_{p}(C')$ (\cite[0.4.8]{A}) 
which naturally 
defines a densely defined operator $T:C_{p}(C') \to C_{p}([0,1])$ 
with domain 
$R((\varphi')^{\sharp}):= (\varphi')^{\sharp}(C_{p}([0,1]))$ 
onto $C_{p}([0,1])$.

\vspace{5mm}

\noindent
Proof of Proposition 2.1.
Our proof is a modification of \cite[Theorem 3.4]{LLP} and is based 
on an analysis of the dual spaces.  We first recall basics of the 
dual space notions \cite[Chap.0]{A}.  
For a Tychonoff space $X$, $L_{p}(X)$ denotes 
the dual space, that is, the space of all continuous 
linear functionals on 
$C_{p}(X)$ endowed with the pointwise convergence topology.  The 
space $L_{p}(X)$ is linearly homeomorphic to the space
\[
\{ \sum_{i=1}^{n}\alpha_{i}x_{i}~|~ \alpha_{i} \in \mathbb{R}, x_{i} \in X, n \in \mathbb{Z}_{\geq 0} \},
\]
where the topology is described below.  Each non-zero point $x$ of $L_{p}(X)$ is uniquely written as
\begin{equation}
x = \sum_{i=1}^{n}\alpha_{i}x_{i}
\end{equation}
where $\alpha_{i} \in \mathbb{R}^{\ast} := \mathbb{R} \setminus \{0\}$ 
and $x_{1}, \ldots, x_{n}$ are mutually distinct points of $X$. The 
number $n$ above is called the length of $x$ and is denoted by 
$\ell(x)$. Let $A_{n}(X) = \{ x \in L_{p}(X)~|~ \ell(x) = n \}$ and 
$B_{n}(X) = \{ x \in L_{p}(X)~|~ \ell(x) \leq n \}$.  
By the definition we have $A_{n}(X) = B_{n}(X) \setminus B_{n-1}(X)$ 
and also we have the equality 
\begin{equation}
L_{p}(X) = \cup_{n=1}^{\infty}A_{n}(X).
\end{equation} 
It is known that 
each $B_{n}(X)$ is closed in $L_{p}(X)$ and 
each point $x = \sum_{i=1}^{n}\alpha_{i}x_{i}$ as in (2.1) has 
a neighborhood basis consisting of the sets
\[
\sum_{i=1}^{n}O_{i}U_{i}
\]
where $O_{i}$'s are open neighborhoods of $\alpha_{i}$'s in $\mathbb{R}^{\ast}$ and $U_{i}$'s are open neighborhoods of $x_i$'s 
respectively 
such that $U_{i} \cap U_{j} = \emptyset$ whenever $i \neq j$. 
It follows from this that 
the map
\begin{equation}
\sigma_{n}:(\mathbb{R}^{\ast})^{n} \times (X^{n}\setminus \Delta_{n}) \to A_{n}(X), ~~((\alpha_{i}),(x_{i})) \mapsto \sum_{i=1}^{n}\alpha_{i}x_{i}, 
\end{equation}
where
$\Delta_{n} = 
\{ (x_{1}, \ldots,x_{n})~|~ x_{i} = x_{j} \mbox{ for some}~i \neq j \}$, 
is a homeomorphism.

Starting the proof of Proposition 2.1, let 
$\xi:E \to C_{p}(X)$ and $\eta:F \to C_{p}(Y)$ be the dense 
inclusions and 
consider the dual maps $\xi^{\ast}:L_{p}(X) \to E^{\ast}$ and 
$\eta^{\ast}:L_{p}(Y) \to F^{\ast}$.  These maps $\xi^{\ast}$ and 
$\eta^{\ast}$ are continuous bijections 
because of the denseness of $E$ and $F$.  It follows from this that 
for each compact set $K$ of 
$L_{p}(X)$, the restriction $\xi^{\ast}|K:K \to \xi^{\ast}(K)$ is a 
homeomorphism.  The same holds for compact sets of $L_{p}(Y)$.
The subset $\Delta_{n}$ above is closed and thus $G_{\delta}$ 
in a metrizable compact space $X^n$.  Then by (2.3), we see that
$A_{n}(X)$ is $\sigma$-compact and is represented as 
the union $A_{n}(X) = \cup_{i=1}^{\infty} A_{n,i}$, 
where each $A_{n,i}$ is homeomorphic to a compact subset of 
$(\mathbb{R}^{\ast})^{n} \times X^{n}$. 
By the above remark, we have
\[
\xi^{\ast}(A_{n}(X)) = \cup _{i =1}^{\infty}\xi^{\ast}(A_{n,i})
\]
and each $\xi^{\ast}(A_{n,i})$ is homeomorphic to $A_{n,i}$.

The composition $T^{\ast} \circ \eta^{\ast}:L_{p}(Y) \to E^{\ast}$ embeds 
$Y$ into $E^{\ast}$ and hence $Y$ is homeomorphic to the subspace
\[
T^{\ast} \eta^{\ast}(Y) \subset \cup_{n,i =1}^{\infty} \xi^{\ast}(A_{n,i}).
\]
Thus we have $Y = \cup_{k = 1}^{\infty}Y_{k}$ where $Y_{k}$ is a compact set such that 
the restriction  

\[
(\xi^{\ast}|A_{n(k),i(k)})^{-1} \circ T^{\ast} \circ \eta^{\ast}|Y_{k}:
Y_{k} \to Y'_{k}
\]
is a homeomorphism of $Y_{k}$ onto a compact subset $Y'_{k}$ of  
$A_{n(k),i(k)}$ for some $n(k)$ and $i(k)$.
Below we show that $\dim Y_{k} = 0$ from which the desired 
conclusion follows by the countable sum theorem 
\cite[Theorem 3.1.8]{En}.

Let $\tilde{Y}_{k} = \sigma_{n(k)}^{-1}(Y_{k}')$ where 
$\sigma_{n(k)}$ is the homeomorphism of (2.3). 
Let $p:(\mathbb{R}^{\ast})^{n(k)} \times X^{n(k)} \to X^{n(k)}$ be the 
projection and consider the restriction $p_{k} := p|\tilde{Y}_{k}$.
As a map defined on the compact set $\tilde{Y}_{k}$, the map 
$p_{k}$ is a closed map onto a zero-dimensional subspace 
$p_{k}(\tilde{Y}_{k})$ of $X^{n(k)}$.
We show that $p_{k}$ is at-most-$n(k)$-to-1 map, that is,
\begin{equation}
|p_{k}^{-1}(x)| \leq n(k)
\end{equation}
for each $x \in X^{n}$. 
Once (2.4) is verified, we obtain the equality $\dim \tilde{Y}_{k} = 0$ 
from the zero-dimensionality of 
$p_{k}(\tilde{Y}_{k})$ by the dimension lowering theorem 
\cite[Theorem 1.12.4]{En}, and therefore we conclude $\dim Y_{k} = 0$, 
as desired.

In what follows, $n(k)$ and $p_{k}$ are simply denoted by $n$ and $p$ respectively.  
In order to verify (2.4), take a point 
$x = (x_{1}, \ldots, x_{n})$ of $X^{n}$ and 
suppose on the contrary that the fiber $p^{-1}((x_{1}, \ldots, x_{n}))$ contains 
$(n+1)$ points $y^{1}, \ldots, y^{n+1}$ of $\tilde{Y}_{k}$. 
For each $j = 1, \ldots, n+1$, 
we may find $\lambda_{ij} \in \mathbb{R}^{\ast}$ such that 
\begin{equation}
y^{j} = \sum_{i=1}^{n}\lambda_{ij}x_{i}.
\end{equation}
By the definition of $\tilde{Y}_{k}$, the definitions of $\sigma_{n(k)}$ 
and the duality, the above (2.5) is rephrased as follows: 
for each $f \in E$ we have
\begin{equation}
(Tf)(y^{j}) = \sum_{i=1}^{n}\lambda_{ij}f(x_{i}), ~~j=1, \ldots,n+1.
\end{equation}
Consider the matrix $\Lambda = (\lambda_{ij})_{1 \leq i,j \leq n}$ of 
size $n \times n$.  The denseness of $F$ with (2.6) implies that the linear map
\[
{\bf v} \mapsto \Lambda {\bf v}, ~~\mathbb{R}^{n} \to \mathbb{R}^{n}
\]
has the dense image.  This implies that $\det \Lambda \neq 0$.  In 
particular we may find $a_{ij} \in \mathbb R$ such that 
\[
(Tf)(y^{n+1}) = \sum_{i=1}^{n}a_{ij}Tf(y^{j}).
\]
Then the image of the evaluation map
\[
e:F \to \mathbb{R}^{n+1}; ~~e(g) = (g(y^{j}))_{1 \leq j \leq n+1}
\]
is contained in an $n$-dimensional subspace of $\mathbb{R}^{n+1}$. 
However since $F$ is dense in $C_{p}(Y)$, the set $e(F)$ must be 
dense in $\mathbb{R}^{n+1}$.  This contradiction 
finishes the proof of (2.4) and thus finishes the proof of 
Proposition 2.1.

\vspace{5mm}

Here we recall some basics on inverse spectra from 
\cite[Chap.1]{Chigogidze}. 
Let  
$\mathcal{S}_{X} = \{X_{\alpha}, p_{\alpha \beta};~A\}$ 
be an inverse system of topologial spaces $X_{\alpha}$,
indexed by a directed set $A$ with the limit space 
$X = \lim_{\leftarrow}\mathcal{S}_{X}$. 
The canonical projection of $X$ to $X_{\alpha}$ is 
denoted by $p_{\alpha}:X \to X_{\alpha}$.
An inverse system 
$\mathcal{S}_{X} = \{X_{\alpha}, p_{\alpha \beta};~A\}$ 
is called a 
{\it factorizing} $\omega$-{\it spectrum} if 
\begin{itemize}
\item[(O1)] each countable chain $C$ of $A$ has the supremum $\sup C \in A$,
\item[(O2)] for each countable chain $B$ of $A$ with $\beta = \sup B$, 
the canonical map 
$\triangle _{\alpha \in B} p_{\alpha, \beta}:X_{\beta} 
\to \lim_{\leftarrow}\{X_{\alpha}, p_{\alpha_{1} \alpha_{2}};~B\}$ 
is a topological embedding, and
\item[(F)] each continuous function $f:X = \lim_{\leftarrow}\mathcal{S}_{X} \to \mathbb{R}$ 
admits an $\alpha \in A$ and $f_{\alpha}:X_{\alpha} \to \mathbb{R}$ such that 
$f = f_{\alpha} \circ p_{\alpha}$. 
\end{itemize}

For a compact Hausdorff space $X$, there exists a factorizing 
$\omega$-spectrum 
$\mathcal{S}_{X} = \{X_{\alpha}, p_{\alpha \beta};~A\}$ with 
$|A| \leq w(X)$
such that $X = \lim_{\leftarrow}\mathcal{S}_{X}$   
and 
\begin{itemize}
\item[(C1)] each $X_{\alpha}$ is a compact metrizable space,
\item[(C2)] each limit projection $p_{\alpha}$ as well as each bonding map 
$p_{\alpha \beta}$ is surjective.
\item[(C3)] the canonical map $\triangle _{\alpha \in B} p_{\alpha \beta}$ 
of (O2) is a surjective homeomorphism.
\end{itemize}
If $\dim X \leq n$, then we may choose above spectrum so that 
$\dim X_{\alpha} \leq n$ for each $\alpha \in A$ 
\cite[Propositions 1.3.5, 1.3.2, and 1.3.10]{Chigogidze}.
For a compact metrizable space $X$, the above spectrum is reduced to 
the trivial system $\{ X, \operatorname{id}_{X} \}$.

\begin{prop}
Let $X = \lim_{\leftarrow} \mathcal{S}_{X}$ and 
$Y = \lim_{\leftarrow}\mathcal{S}_{Y}$ be 
compact Hausdorff spaces  
which are the limits of 
factorizing $\omega$-spectra 
$\mathcal{S}_{X} = \{ X_{\alpha}, p_{\alpha_{1} \alpha_{2}};~A \}$ and 
$\mathcal{S}_{Y} = \{ Y_{\alpha}, p_{\beta_{1} \beta_{2}};~B \}$ 
satisfying the conditions (C1)-(C3) with the projections 
$p_{\alpha}:X \to X_{\alpha}$ and $q_{\beta}:Y \to Y_{\beta}$ 
respectively.

Let $T:C_{p}(X) \to C_{p}(Y)$ be a linear continuous opeartor.
\begin{itemize}
\item[(1)] For each $\alpha$, there exist $\beta = \beta(\alpha)$ 
and a densely defined operator $T_{\alpha,\beta}:C_{p}(X_{\alpha}) \to C_{p}(Y_{\beta})$ such that
\[
T \circ p_{\alpha}^{\sharp} = q_{\beta}^{\sharp} \circ T_{\alpha, \beta}.
\]
For each $\beta_{0} \in B$, we may choose the above $\beta(\alpha)$ 
so that $\beta (\alpha) \geq \beta_{0}$. 
\item[(2)] If moreover $T$ is surjective, then for each 
$\alpha_{0} \in A$ and for each $\beta_{0} \in B$, we may choose 
$T_{\alpha,\beta}$ so that $\alpha \geq \alpha_{0}, \beta \geq \beta_{0}$ and $T_{\alpha,\beta}$ has a dense image.
\end{itemize}
\end{prop}

The following diagram illustrates the operator $T_{\alpha, \beta}$.

$$
\begin{CD}
C_{p}(X) @>T>> C_{p}(Y) \\
@Ap_{\alpha}^{\sharp}AA @ AAq_{\beta}^{\sharp}A\\
C_{p}(X_{\alpha}) @>T_{\alpha,\beta}>> C_{p}(Y_{\beta})
\end{CD}
$$

Proof.
Let $f_{\alpha} \in C_{p}(X_{\alpha})$ and consider the composition 
$f_{\alpha} \circ p_{\alpha}$ whose image by $T$ is factorized as
\begin{equation}
T(f_{\alpha} \circ p_{\alpha}) = g_{\beta} \circ q_{\beta}
\end{equation}
for some $\beta$ and $g_{\beta} \in C_{p}(Y_{\beta})$.  Observe that 
we may choose the above $\beta$ as large as we wish.  An important 
observation here is that, because $q_{\beta}$ is surjective, 
\begin{equation}
g_{\beta} ~\mbox{ is uniquely determined by} ~f_{\alpha} ~\mbox{ and } ~\beta.
\end{equation}

Fix an arbitrary $\alpha$ and notice that $C_{p}(X_{\alpha})$ is 
separable (\cite[Theorem I.1.5]{A}).  Take a countable dense 
set $D = \{ f_{\alpha,i} \}$ of $C_{p}(X_{\alpha})$. For each $i$
take a $\beta_{i}$ and $g_{i} \in C_{p}(Y_{\beta_{i}})$ such that 
\[
T(f_{\alpha,i} \circ p_{\alpha}) = g_{i} \circ q_{\beta_{i}}.
\]
Let $\beta = \sup_{i}\beta_{i} \in B$.  By the $\omega$-continuity of 
the spectrum (C3) we have
\[
Y_{\beta} = \lim_{\leftarrow} \{ Y_{\beta_{i}}, q_{\beta_{i} \beta_{j}}\}
\]
with the projection $q_{\beta_{i},\beta}:Y_{\beta} \to Y_{\beta_i}$.
For each $i$, the function 
$g_{\beta,i} = q_{\beta _{i} \beta}^{\sharp}(g_{i}) \in C_{p}(Y_{\beta})$ satisfies
 \[
T(f_{\alpha,i} \circ p_{\alpha}) = g_{\beta,i} \circ q_{\beta}.
\]
Define 
$T_{\alpha,\beta}:D \to C_{p}(Y_{\beta})$ by
\[
T_{\alpha,\beta}(f_{\alpha,i}) = g_{\beta,i} 
\]
and let $D(T_{\alpha,\beta}) = \mathrm{span}D$ which is a dense 
subspace of $C_{p}(X_{\alpha})$. 
We make use of the uniqueness (2.8)  
to extend $T_{\alpha,\beta}$ to a linear 
operator $T_{\alpha,\beta}:D(T_{\alpha, \beta}) \to C_{p}(Y_{\beta})$ 
as follows.

For $f = \sum_{i} \lambda_{i}f_{\alpha, i} \in D(T_{\alpha, \beta})$, 
all but finitely many $\lambda_{i}$'s being zero,  
let $T_{\alpha, \beta}(f) = \sum_{i}\lambda_{i} g_{\beta,i}$.
Then $T_{\alpha, \beta}$ is well-defined and a linear map. 
Indeed, suppose that $\sum_{i} \lambda_{i}f_{\alpha, i} = 
\sum_{i} \mu_{i}f_{\alpha, i}$.  Then we have the euality 
$\sum_{i} \lambda_{i}f_{\alpha, i} \circ p_{\alpha} = 
\sum_{i} \mu_{i}f_{\alpha, i} \circ p_{\alpha}$ and hence we have 
that 
$T(\sum_{i} \lambda_{i}f_{\alpha, i} \circ p_{\alpha}) = 
T(\sum_{i} \mu_{i}f_{\alpha, i} \circ p_{\alpha})$.
By linearity of $T$ we see that 
\begin{eqnarray*}
\sum_{i} \lambda_{i} g_{\beta, i} \circ q_{\beta} &=&  
\sum_{i} \lambda_{i} T(f_{\alpha, i} \circ p_{\alpha}) 
= T(\sum_{i} \lambda_{i} f_{\alpha, i} \circ p_{\alpha}) \\
&=& T(\sum_{i} \mu_{i} f_{\alpha, i} \circ p_{\alpha})
= \sum_{i} \mu_{i} T(f_{\alpha, i} \circ p_{\alpha})\\
&=& \sum_{i} \mu_{i} g_{\beta, i} \circ q_{\beta},
\end{eqnarray*}
from which we see that 
$\sum_{i} \lambda_{i} g_{\beta, i} = 
\sum_{i} \mu_{i} g_{\beta, i}$ 
by the surjectivity of $q_{\beta}$.
This proves that $T_{\alpha, \beta}$ is well-defined. 
The same argument is applied to prove that $T_{\alpha, \beta}$ 
is a linear map.

Finally we verify the continuity of $T_{\alpha,\beta}$.
For a finite set $F$ of a space $Z$, for an $\epsilon > 0$ 
and for a function 
$h \in C_{p}(Z)$, let 
\[
<h, F, \epsilon> = \{ u \in C(Z)~|~ |u(p) - h(p) | < \epsilon  ~~\mbox{for each} ~~p \in F\}.
\]
Fix an $\epsilon >0$ and a finite subset $F_{\beta}$ of $Y_{\beta}$. 
For each 
$y_{\beta} \in F_{\beta}$, take $y \in Y$ such that 
$q_{\beta}(y) = y_{\beta}$ by 
the surjectivity of $q_{\beta}$.
Let $F$ be the resulting finite subset of $Y$.  By the continuity of $T$, 
we may find a finite subset $E$ of $X$ and $\delta>0$ such that
\[
T(<f_{\alpha} \circ p_{\alpha}, E, \delta>) \subset 
<T(f_{\alpha} \circ p_{\alpha}), F, \epsilon>.
\]
Let $E_{\alpha} = p_{\alpha}(E)$.  We verify the inclusion
\[
T_{\alpha,\beta}(<f_{\alpha}, E_{\alpha}, \delta> \cap D(T_{\alpha,\beta})) \subset <T_{\alpha,\beta}(f_{\alpha}), F_{\beta}, \epsilon>.
\]
Indeed for $f'_{\alpha} \in <f_{\alpha}, E_{\alpha}, \delta>$, we have 
$|f'_{\alpha}(p_{\alpha}(x)) - f_{\alpha}(p_{\alpha}(x))| < \delta$ 
for each $x \in E$ and 
hence $f'_{\alpha} \circ p_{\alpha} \in <f_{\alpha} \circ p_{\alpha}, E, \delta>$. For each $y \in F$ we have 
\[
|T(f'_{\alpha} \circ p_{\alpha})(y) - T(f_{\alpha} \circ p_{\alpha})(y)| 
< \epsilon, 
\]
which implies 
\[
|T_{\alpha,\beta}(f'_{\alpha})(y_{\beta}) - T_{\alpha,\beta}(f_{\alpha})(y_{\beta})| < \epsilon
\]
for each $y_{\beta} \in F_{\beta}$ 
by the equality 
$T_{\alpha \beta}(f_{\alpha}) \circ q_{\beta} = T(f_{\alpha}\circ p_{\alpha})$. 
Thus we see that 
$T_{\alpha,\beta}(f'_{\alpha}) \in <T_{\alpha,\beta}(f_{\alpha}), F_{\beta}, \epsilon>$.  
This proves the continuity of $T_{\alpha, \beta}$ and hence completes the proof of (1).

(2) Now assume that $T$ is surjective. For an arbitrary 
$\alpha$ with $\alpha \geq \alpha_{0}$, 
take $\beta = \beta(\alpha) \geq \beta_{0}$ and 
$T_{\alpha,\beta}:C_{p}(X_{\alpha}) \to C_{p}(Y_{\beta})$ as 
in (1).  The operator $T_{\alpha,\beta}$ is defined on a dense subspace 
$D(T_{\alpha, \beta})$ of $C_{p}(X_{\alpha})$ and satisfies 
\begin{equation}
T \circ p_{\alpha}^{\sharp} = q_{\beta}^{\sharp} \circ T_{\alpha,\beta}.
\end{equation}
Take a countable dense subset $E = \{ g_{i} \}$ of $C_{p}(Y_{\beta})$. 
Since $T$ is surjective, for each $i$ there exists $f_{i} \in C_{p}(X)$ such that 
$T(f_{i}) = g_{i} \circ q_{\beta}$. The function $f_i$ factorizes as 
\[
f_{i} = f_{\alpha_{i}} \circ p_{\alpha_{i}}
\]
for some $\alpha_{i} \in A$. 
Let $\alpha (1) = \sup_{i} \alpha_{i}$. The function 
$f_{\alpha (1),i} = f_{\alpha _{i}} \circ p_{\alpha_{i} \alpha (1)}$ satisfies 
$f_{i} = f_{\alpha (1),i} \circ p_{\alpha (1)}$ and hence 
\begin{equation}
T(f_{\alpha (1),i} \circ p_{\alpha (1)}) = g_{i} \circ q_{\beta}.
\end{equation}
Choose a countable dense subset $D_{1}$ of $C_{p}(X_{\alpha (1)})$ such 
that 
$D_{1} \supset p_{\alpha, \alpha(1)}^{\sharp}(D) \cup \{f_{\alpha (1),i} \}$.
Repeat the procedure of (1) with the use of $D_{1}$ to
find $\beta (1) > \beta$ and 
a linear map $T_{\alpha (1),\beta (1)}:C_{p}(X_{\alpha (1)}) \to C_{p}(Y_{\beta (1)})$ 
densely defined on the linear span $D(T_{\alpha (1),\beta(1)})$ of $D_{1}$ such that 
\begin{equation}
T(f_{\alpha(1)} \circ p_{\alpha(1)}) = T_{\alpha(1),\beta(1)} \circ q_{\beta(1)}
\end{equation}
for each $f_{\alpha(1)} \in D(T_{\alpha(1),\beta(1)})$. By (2.10), (2.11), and the surjectivity 
of $q_{\beta(1)}$ we see that $T_{\alpha(1),\beta(1)}(f_{\alpha(1),i}) = g_{i}$ 
for each $i$.  Thus we have the inclusion
$R(T_{\alpha(1),\beta(1)}) \supset 
q^{\sharp}_{\beta, \beta(1)}(E).$

Let us summarize the above properties of $T_{\alpha(1),\beta(1)}$:
\begin{itemize}
\item[(1.1)] $D(T_{\alpha (1),\beta(1)}) \supset 
p^{\sharp}_{\alpha,\alpha(1)}(D(T_{\alpha,\beta}))$ and
\item[(1.2)] $R(T_{\alpha(1),\beta(1)}) \supset 
q^{\sharp}_{\beta,\beta(1)}(E)$.
\end{itemize}

We continue this process to obtain increasing sequences of 
indices 
$\alpha(1) < \cdots < \alpha(n) < \cdots$, 
$\beta(1) <  \cdots < \beta(n) < \cdots $,
countable dense subsets $D_{n}$ of $C_{p}(X_{\alpha(n)})$, 
$E_{n}$ of $C_{p}(Y_{\beta(n)})$,
 and 
a sequence of linear continuous operators 
$\{ T_{\alpha(n),\beta(n)}:C_{p}(X_{\alpha(n)}) \to C_{p}(Y_{\beta(n)})\}$, 
each $T_{\alpha(n),\beta(n)}$ being defined on $D(T_{\alpha(n),\beta(n)})$, such that
\begin{itemize}
\item[(n.1)] $D(T_{\alpha(n+1),\beta(n+1)}) \supset p_{\alpha(n)\alpha(n+1)}^{\sharp}(D(T_{\alpha(n),\beta(n)}))$,
\item[(n.2)] $R(T_{\alpha(n+1),\beta(n+1)}) \supset 
q_{\beta(n)\beta(n+1)}^{\sharp}(E_{n})$, and
\item[(n.3)] $T_{\alpha(n+1),\beta(n+1)} \circ 
p_{\alpha(n) \alpha(n+1)}^{\sharp} = q_{\beta(n) \beta(n+1)}^{\sharp} \circ T_{\alpha(n),\beta(n)}$ 
on $D(T_{\alpha(n),\beta(n)})$.
\end{itemize}
Let $\alpha_{\infty} = \sup_{n} \alpha(n)$ and 
$\beta_{\infty} = \sup_{n} \beta(n)$ and let 
$D_{\infty} = 
\cup_{n} 
p_{\alpha(n) \alpha_{\infty}}^{\sharp}(D(T_{\alpha(n),\beta(n)}))$.
Then 
\begin{equation}
D_{\infty}~\mbox{ is dense in } ~C_{p}(X_{\alpha_{\infty}}).
\end{equation}
Indeed $X_{\alpha_{\infty}} = \lim_{\leftarrow}X_{\alpha(n)}$ by the 
$\omega$-continuity (C3).  Hence for each $f \in C_{p}(X_{\alpha_{\infty}})$ 
and for each $\epsilon > 0$, there exist $\alpha(n)$ and 
$f_{\alpha(n)} \in C_{p}(X_{\alpha(n)})$ such that 
\[
\sup_{x_{\infty} \in X_{\alpha_{\infty}}}
 |f(x_{\infty}) - f_{\alpha(n)} (p_{\alpha(n)}(x_{\infty}))| < \epsilon.
\]
The function $f_{\alpha(n)}$ is approximated arbitrarily closely 
by functions of 
$D(T_{\alpha(n),\beta(n)})$ with respect to the pointwise convergence topology on
$C_{p}(X_{\alpha(n)})$.  Hence the function $f$ is approximated arbitrarily closely 
by functions from the set 
$\cup _{n} p_{\alpha(n)\alpha_{\infty}}^{\sharp}(D(T_{\alpha(n),\beta(n)})) 
= D_{\infty}$.

An operator $T_{\alpha_{\infty},\beta_{\infty}}$ 
is defined by
\[
T_{\alpha_{\infty},\beta_{\infty}}(p_{\alpha(n), \alpha_{\infty}}^{\sharp}(f_{\alpha(n)})) = 
T_{\alpha(n),\beta(n)}(f_{\alpha(n)}) \circ q_{\beta(n),\beta_{\infty}}
\]
for $f_{\alpha(n)} \in D(T_{\alpha(n), \beta(n)})$.  
The equality (n.3) guarantees that $T_{\alpha_{\infty},\beta_{\infty}}$ 
is well defined. Now we have the inclusion
$T_{\alpha_{\infty},\beta_{\infty}}(
p^{\sharp}_{\alpha_{\infty} \alpha(n)}(D(T_{\alpha(n),\beta(n)}))) 
\supset q_{\beta(n) \beta_{\infty}}(E_{n})$.  
Hence $R(T_{\alpha_{\infty},\beta_{\infty}}) \supset E_{\infty}:= 
\cup_{n}q^{\sharp}_{\beta(n)\beta_{\infty}}(E_{n})$.  
Using the same argument as in (2.12) we see that $E_{\infty}$ 
is dense in $C_{p}(Y_{\beta_{\infty}})$.

This proves (2).

\vspace{5mm}
\noindent
Proof of Theorem 1.1.  Let $T:C_{p}(X) \to C_{p}(Y)$ be a linear continuous surjection 
where $X$ and $Y$ are compact Hausdorff spaces.  Let 
$\mathcal{S}_{X} = \{ X_{\alpha}, p_{\alpha_{1} \alpha_{2}}; ~A \}$ and 
$\mathcal{S}_{Y} = \{ Y_{\alpha}, p_{\beta_{1} \beta_{2}}; ~B \}$
be factorizing $\omega$-spectra such that $X = \lim_{\leftarrow}\mathcal{S}_{X}$ and 
$Y= \lim_{\leftarrow}\mathcal{S}_{Y}$.  By the assumption $\dim X = 0$, 
we may assume that $\dim X_{\alpha} = 0$ for each $\alpha \in A$.  
For each $\beta \in B$, we
apply Proposition 2.2 to find a densely defined 
operator 
$T_{\alpha,\beta(\alpha)}:C_{p}(X_{\alpha}) 
\to C_{p}(Y_{\beta(\alpha)})$ 
such that $\beta(\alpha) \geq \beta$ and 
$R(T_{\alpha,\beta(\alpha)})$ is dense in 
$C_{p}(Y_{\beta(\alpha)})$. By Propostion 2.1, we have 
$\dim Y_{\beta(\alpha)} = 0$ and hence 
the set $\{ \beta \in B ~|~\dim Y_{\beta} = 0 \}$ forms 
a cofinal subset of $B$. We therefore obtain 
$\dim Y = 0$.  
This completes the proof of Theorem 1.1.

\begin{rem}
Theorem 1.1 holds under a weaker hypothesis that $X$ is a zero-dimensional 
pseudocompact Tychonoff space and $Y$ is a compact Hausdorff space.  
\end{rem}
\noindent
To see the above, assume that $X$ is such a space and $T:C_{p}(X) \to C_{p}(Y)$ 
is a linear continuous surjection onto $C_{p}(Y)$ where $Y$ is compact.
The inclusion 
$h:X \to \beta X$ of $X$ into the 
Stone-{\v C}ech compactification $\beta X$ of $X$ 
induces a linear continuous 
surjection  
$h^{\sharp}:C_{p}(\beta X) \to C_{p}(X)$.  Since 
$\dim \beta X = 0$, we may apply Theorem 1.1 to 
the composition $T \circ h^{\sharp}$ to conclude 
that $\dim Y = 0$.

\vspace{5mm}
\noindent
Proof of Proposition 1.2.  Let $X$ be a Bing compactum, let 
$Y$ be a hereditarly locally connected compact metrizable space, 
and let $T:C_{p}(X) \to C_{p}(Y)$ be a densely defined 
linear continuous operator with a dense image. 
First recall that each continuous map $\varphi:Z \to B$ 
of a locally connected continuum $Z$ to a Bing compactum $B$ must be 
a constant map.  Proceeding as in the proof of Proposition 2.1, 
we see that 
the space $Y$ is the countable union
$Y = \cup_{i =1}^{\infty} Y_{i}$ 
such that each $Y_{i}$ is homeomorphic to a subspace $\tilde{Y}_{i}$ of 
$(\mathbb{R}^{\ast})^{n} \times X^{n}$, and such that the projection 
$p_{i}:\tilde{Y}_{i} \to X^{n}$ is a finite-to-one map.
     
By the assumption, each component $C$ of $Y_{i}$ is locally connected.  
Applying the above remark to the composition of the map $p_{i}|C:C \to X^{n}$ 
with the projection $X^{n} \to X$ onto any factor, we see that 
$p_{i}|C$ must be a constant map.  This implies that $C$ is contained in 
a fiber of $p_{i}$, hence $C$ is a finite set and thus $C$ is a singleton.  This 
proves that $Y_{i}$ is totally disconnected which is equivalent 
to $\dim Y_{i} = 0$ by the compact metrizability of $Y_{i}$. 
By the countable sum theorem \cite[Theorem 3.1.8]{En}, 
we conclude that $\dim Y = 0$. 

\section{Remarks and Problems}

The proof of Propostion 2.1 relies on both the compactness and the 
metrizability of the spaces involved, while Remark 2.4 naturally 
raises the following problem.

\begin{prob}
Does Theorem 1.1 hold for Tychonoff spaces $X$ and $Y$?
\end{prob}

Recently Krupski and Marciszewski \cite{KM} proved that there exists a 
metrizable space $X$ such that $C_{p}(X)$ is not homeomorphic to 
$C_{p}(X) \times C_{p}(X) \simeq C_{p}(X \oplus X)$, negatively 
answering 
a long standing problem of Arkhangel'skii.  In the paper above 
they also showed that there exists no linear continuous 
surjection $C_{p}(M) \to C_{p}(M) \times C_{p}(M)$ for 
the Cook continuum $M$. Recall that the Cook continuum $M$ 
is a hereditarily indecomposable continuum 
such that, for each non-degenerate subcontinuum $C$ of $M$, 
every continuous map $f:C \to M$ is either the identity map  $\operatorname{id}$ 
or a constant map \cite{C}. 
Since $M \oplus M$ is a subspace of $M \times M$, we see 
that there exists no linear 
continuous surjection $C_{p}(M) \to C_{p}(M \times M)$. 
On the other hand as has been pointed out by 
Marciszewski (a private communication), we have the following. 

\begin{rem}
For the pseudo-arc $P$, there exists a topological linear isomorphism 
$C_{p}(P) \simeq C_{p}(P \oplus P)$. 
\end{rem} 

Here we give a sketch of the above.  Take a non-degenerate subcontinuum 
$Q$ of $P$ and let $P/Q$ be the quotient space obtained from $P$ by 
shrinking $Q$ into a point. The space $P/Q$ is a metrizable continuum 
as well. The projection $P \to P/Q$ is a 
monotone map, and monotone maps preserve the hereditary 
indecomposability and the arc-likeness (that is, being represented 
by the limit of an inverse sequence of [0,1]), from the 
characterization of the pseudo-arc, we see that $P/Q$ is homeomorphic to $P$.
Then we obtain, by \cite[Corollary 6.6.13]{vM}, the following linear topological isomorphisms
\[
C_{p}(P) \simeq C_{p}(P/Q) \times C_{p}(Q) \simeq C_{p}(P \oplus P).
\]
This proves the desired result.

The following problem naturally arises.

\begin{prob}
Does there exist a linear continuous surjection 
$C_{p}(P) \to C_{p}(P \times P)$?
\end{prob} 

Here we notice that, if a finite-dimensional 
compact metrizable space $X$ contains a homeomorphic copy of $[0,1]$, 
then there does exist a linear continuous surjection 
$C_{p}(X) \to C_{p}(X \times X)$.  
In order to see this, first observe that an arbitrary embedding 
$h:[0,1] \to X$ induces a linear continuous 
surjection $h^{\sharp}:C_{p}(X) \to C_{p}([0,1])$.  Since 
$X$ is finite-dimensional, compact and metrizable, 
there exists a linear continuous surjection 
$T:C_{p}([0,1]) \to C_{p}(X)$ \cite{LMP}.  Then the  
composition $T \circ h^{\sharp}$ is the desired 
linear continuous surjection.

\begin{rem}
Theorem 1.1 does not hold when we drop the assumption of linearity of 
the operator $T$.  
\end{rem}
In fact we can prove the following:
Let $X$ be a compact metrizable space and let $S$ be the convergent 
sequence, that is, the space homeomorphic to 
$\{0\} \cup \{\frac{1}{n}~|~n \in \omega \}$.  Then there exists a continuous surjection $H:C_{p}(S) \to C_{p}(X)$.

\vspace{5mm}
\noindent
Proof. The argument below is extracted from \cite[Proposition 5.4]{KM}. 
The space $C_{p}(S)$ contains a closed homeomorphic copy of $J$, the space of 
irrationals.  Since the Banach space $C(X)$ of all real-valued 
continuous functions with sup norm is separable, we obtain a 
continuous surjection $J \to C(X)$ which extends to a continuous 
surjection $H:C_{p}(S) \to C(X)$. The map $H$, regarded as a map to 
$C_{p}(X)$, is the desired continuous surjection 
$H:C_{p}(S) \to C_{p}(X)$.

\begin{flushleft}
Kazuhiro Kawamura\\
Institute of Mathematics\\
University of Tsukuba\\
Tsukuba, Ibaraki 305-8071\\
Japan\\
kawamura@math.tsukuba.ac.jp\\
\end{flushleft}
\begin{flushleft}
Arkady Leiderman\\
Department of Mathematics\\
Ben-Gurion University of the Negev\\
Beer-Sheva\\ 
Israel\\
arkady@math.bgu.ac.il\\
\end{flushleft}

\end{document}